\documentclass[12pt]{article}

\usepackage{amsthm}
\usepackage{amsmath}
\usepackage{amssymb}

\usepackage{graphicx}
\usepackage{float}

\newtheorem*{theorem}{Theorem}
\newtheorem{ntheorem}{Theorem}
\newtheorem*{definition}{Definition}
\newtheorem*{proposition}{Proposition}
\newtheorem{nproposition}{Proposition}
\newtheorem*{corollary}{Corollary}

\newcommand{\R}{\mathbb{R}}

\newcommand{\w}{\omega}
\newcommand{\B}{{\cal B}}

\newcommand{\T}{{\cal T}}
\newcommand{\F}{{\cal F}}

\newcommand{\rank}{r}

\begin{document}

\pagestyle{plain}

\title{The Bergman Complex of a Matroid \\ and Phylogenetic Trees}

\author{\begin{tabular}{ccc}
Federico Ardila & & Caroline J. Klivans
\end{tabular}}
\date{}
\maketitle

\maketitle

\begin{abstract}
We study the Bergman complex $\B(M)$ of a matroid $M$: a
polyhedral complex which arises in algebraic geometry, but which
we describe purely combinatorially. We prove that a natural
subdivision of the Bergman complex of $M$ is a geometric
realization of the order complex of its lattice of flats.  In
addition, we show that the Bergman fan $\widetilde{\B}(K_n)$ of
the graphical matroid of the complete graph $K_n$ is homeomorphic
to the space of phylogenetic trees $\T_n$.
\end{abstract}

\section{Introduction}

In~\cite{Bergman}, Bergman defined the \emph{logarithmic
limit-set} of an algebraic variety in order to study its
exponential behavior at infinity. We follow \cite{Sturmfels} in
calling this set the \emph{Bergman complex} of the variety.
Bergman conjectured that this set is a finite, pure polyhedral
complex. He also posed the question of studying the geometric
structure of this set; \emph{e.g.}, its connectedness, homotopy,
homology and cohomology. Bieri and Groves first proved the
conjecture in~\cite{Bieri} using valuation theory.

Recently, Bergman complexes have received considerable attention
in several areas, such as tropical algebraic geometry and
dynamical systems. They are the \emph{non-archimedean amoebas} of
\cite{Einsiedler} and the \emph{tropical varieties} of
\cite{Speyer, Sturmfels}. In particular, Sturmfels
\cite{Sturmfels} gave a new description of the Bergman complex and
an alternative proof of Bergman's conjecture in the context of
Gr\"{o}bner basis theory. Moreover, when the variety is a linear
space, so the defining ideal $I$ is generated by linear forms, he
showed that the Bergman complex can be described solely in terms
of the matroid associated to the linear ideal.

Sturmfels used this description to define the Bergman complex of
an arbitrary matroid, and suggested studying its combinatorial,
geometric and topological properties \cite{Sturmfels}. The goal of
the paper is to undertake this study.

In Section \ref{section:minimum} we study the collection of bases
of minimum weight of a matroid with respect to a weight vector. We
show that this collection is itself the collection of bases of a
matroid, and we give several descriptions of it.

In Section \ref{section:bergman} we prove the main result of the
paper. We show that, appropriately subdivided, the Bergman complex
of a matroid $M$ is the order complex of the proper part of the
lattice of flats $L_M$ of the matroid. These order complexes are
well-understood objects \cite{Bjorner}, and an immediate corollary
of our result is an answer to the questions of Bergman and
Sturmfels about the geometry of $\B(M)$ in this special case. The
Bergman complex of an arbitrary matroid $M$ is a finite, pure
polyhedral complex. In fact, it is homotopy equivalent to a wedge
of $(r-2)$-dimensional spheres, where $r$ is the rank of $M$.

In Section \ref{section:trees}, we take a closer look at the
Bergman complex of the graphical matroid of the complete graph
$K_n$.  We show that the Bergman fan $\widetilde{\B}(K_n)$ is
exactly the space of ultrametrics on $[n]$, which is homeomorphic
to the space of phylogenetic trees of \cite{Billera}. As a
consequence, we show that the order complex of the proper part of
the partition lattice $\Pi_n$ is a subdivision of the link of the
origin of this space. This provides a new explanation and a
strengthening for the known result that these two simplicial
complexes are homotopy equivalent \cite{Robinson, Sundaram,
Trappmann, Vogtmann, Wachs}.

Although we have tried to keep the presentation as self-contained
as possible, some familiarity with the basic notions of matroid
theory will be useful throughout the paper. For the relevant
definitions, we refer the reader to \cite{Oxley}.

\section{The bases of minimum weight of a matroid} \label{section:minimum}

Let $M$ be a matroid of rank $r$ on the ground set
$[n]=\{1,2,\ldots,n\}$, and let $\w \in \R^n$. Regard $\w$ as a
weight function on $M$, so that the weight of a basis $B=\{b_1,
\ldots, b_r\}$ of $M$ is given by $\w_B= \w_{b_1} +\w_{b_2} +
\cdots + \w_{b_r}$.

Let $M_{\w}$ be the collection of bases of $M$ having minimum
$\w$-weight. This is one of the central objects of our study, and
we wish to understand it from three different points of view:
geometric, algorithmic and matroid theoretic.

\medskip

\noindent Geometrically, we can understand $M_{\w}$ in terms of
the matroid polytope. We will use the following characterization
of matroid polytopes, due to Gelfand and Serganova:
\begin{theorem}\cite[Theorem 1.11.1]{Borovik}
Let $S$ be a collection of $r$-subsets of $[n]$. Let $P_S$ be the
polytope in $\R^n$ with vertex set $\{ e_{b_1} + \cdots + e_{b_r}
\, | \, \{b_1, \ldots, b_r\} \in S\}$, where $e_i$ is the $i$-th
unit vector. Then $S$ is the collection
of bases of a matroid if and only if every edge of $P_S$ is a
translate of the vector $e_i - e_j$ for some $i,j \in [n]$.
\end{theorem}

Let $P_M$ be the matroid polytope of $M$. We can now think of $\w$
as a linear functional in $\R^n$. The bases in $M_{\w}$ correspond
to the vertices of $P_M$ which minimize the linear functional
$\w$. Their convex hull is $P_{M_{\w}}$, the face of $P_M$ where
$\w$ is minimized. It follows that the edges of
$P_{M_{\w}}$, being edges of $P_{M}$ also, are parallel to vectors
of the form $e_i-e_j$. Therefore $M_{\w}$ is the collection of
bases of a matroid.

\medskip

\noindent Algorithmically, matroids have the property that
their $\w$-minimum bases are precisely the possible outputs of the
greedy algorithm: Start with $B=\emptyset$. At each
stage, look for an $\w$-minimum element of $[n]$ which can be
added to $B$ without making it dependent, and add it. After
$\rank$ steps, output the basis $B$. \cite[Theorem 1.8.5]{Oxley}
%
%
%
%

\begin{definition}
Given $\omega \in \R^n$, let $\F(\omega)$ denote the unique flag
of subsets
$$
\emptyset=:F_0 \subset F_1 \subset \cdots \subset F_{k} \subset
F_{k+1}:=[n]
$$
for which $\omega$ is constant on each set $F_i - F_{i-1}$ and has
$\omega|_{F_i-F_{i-1}} < \omega|_{F_{i+1}-F_i}$. The \emph{weight
class} of a flag $\F$ is the set of vectors $\w$ such that $\F(\w)
= \F$.

\end{definition}

We can describe weight classes by their defining equalities and
inequalities. For example, one of the weight classes in $\R^5$ is
the set of vectors $\w$ such that $\w_1=\w_4 < \w_2 < \w_3 =
\w_5$. It corresponds to the flag $\{\emptyset \subset \{1,4\}
\subset \{1,2,4\} \subset \{1,2,3,4,5\}\}$.

\begin{nproposition}  \label{prop:Mw1}
If $\omega$ is in the weight class of $\F=\{\emptyset =: F_0
\subset \ldots \subset F_{k+1} := [n]\}$, then the
$\omega$-minimum bases of $M$ are exactly those containing
$\rank(F_i)-\rank(F_{i-1})$ elements of $F_i-F_{i-1}$, for each
$i$. Consequently, $M_{\w}$ depends only on $\F$, and we call it
$M_{\F}$.
\end{nproposition}

\begin{proof}
The greedy algorithm picks $\rank(F_1)$ elements of the lowest
weight, until it reaches a basis of $F_1$; then it picks
$\rank(F_2)-\rank(F_1)$ elements of the second lowest weight,
until it reaches a basis of $F_2$, and so on. Therefore, the
possible outputs of the algorithm are precisely the ones
described.
%
\end{proof}

\medskip

\noindent Matroid theoretically, $M_{\w}$ can be constructed as a
direct sum of minors of $M$, and its lattice of flats $L_{M_{\w}}$
can be constructed from intervals of $L_M$. Let $M|S$ and $M/S$
denote, respectively, the restriction and contraction of the
matroid $M$ to a subset $S$ of its ground set. Then we have:

\begin{nproposition} \label{prop:Mw2}
If $\F=\{\emptyset =: F_0 \subset \ldots \subset F_{k+1} :=
[n]\}$, then
$$
M_{\F}=\bigoplus_{i=1}^{k+1} \, (M|F_i)/F_{i-1} \qquad
\textrm{and} \qquad L_{M_{\F}} \cong \prod_{i=1}^{k+1}
[F_{i-1},F_i].
$$
\end{nproposition}

\begin{proof}
After $\rank(F_{i-1})$ steps, the greedy algorithm has chosen a
basis of $F_{i-1}$. In the following $\rank(F_{i})-\rank(F_{i-1})$
steps, it needs to choose elements which, when added to $F_{i-1}$,
give a basis of $F_i$. The possible choices are, precisely, the
bases of $(M|F_i)/F_{i-1}$. The first equality follows, and the
second one follows from it.
%
\end{proof}

\section{The Bergman complex}\label{section:bergman}

We now define the two main objects of study of this paper.

\begin{definition}
The \emph{Bergman fan} of a matroid $M$ with ground set $[n]$ is
the set
\[
\widetilde{\B}(M) := \{ \w \in \R^n \,\, : \,\, M_{\w} \,\,
\text{has no loops} \}.
\]
The \emph{Bergman complex} of $M$ is
\[
\B(M) := \{ \w \in S^{n-2} \,\, : \,\, M_{\w} \,\, \text{has no
loops} \},
\]
where $S^{n-2}$ is the sphere $ \{\,\w \in \R^n \,\, : \,\, \w_1 +
\cdots + \w_n = 0 \, , \, \w_1^2 + \cdots + \w_n^2 = 1 \}.$
\end{definition}

For the moment, we are slightly abusing notation by calling these
two objects a \emph{fan} and a \emph{complex}. We will very soon
see that they are a polyhedral fan and a spherical polyhedral
complex, respectively; this justifies their name. We will
concentrate on the Bergman complex, but the same arguments apply
to the Bergman fan.

Since the matroid $M_{\w}$ only depends on the weight class that
$\w$ is in, the Bergman complex of $M$ is a disjoint union of the
following weight classes:

\begin{definition}
The weight class of a flag $\F$ is \emph{valid} for $M$ if
$M_{\F}$ has no loops.
\end{definition}

We will study two polyhedral subdivisions of $\B(M)$, one of which
is clearly finer than the other.

\begin{definition}
The \emph{fine subdivision} of $\B(M)$ is the subdivision of
$\B(M)$ into valid weight classes: two vectors $u$ and $v$ of
$\B(M)$ are in the same class if and only if $\F(u) = \F(v)$.

The \emph{coarse subdivision} of $\B(M)$ is the subdivision of
$\B(M)$ into $M_{\w}$-equivalence classes: two vectors $u$ and $v$
of $\B(M)$ are in the same class if and only if $M_u = M_v$.
\end{definition}

Recall that the \emph{order complex} $\Delta(P)$ of a poset $P$ is
the simplicial complex whose vertices are the elements of $P$, and
whose faces are the chains of $P$.

\begin{ntheorem}\label{theorem:main}
The weight class of a flag $\F$ is \emph{valid} for $M$ if and
only if $\F$ is a flag of flats of $M$.
Therefore, the fine subdivision of the Bergman complex $\B(M)$ is
a geometric realization of $\Delta(L_M-\{\,\hat{0}\,,
\hat{1}\,\}\,)$, the order complex of the proper part of the
lattice of flats of $M$.

\end{ntheorem}

\begin{proof}
Assume $F_i$ in $\F$ \emph{is not} a flat of $M$, so there exists
some $e \notin F_i$ in the closure $\overline{F_i}$. By
Proposition \ref{prop:Mw1}, any basis $B$ in $M_{\F}$ contains
$\rank(F_i)$ elements of $F_i$; since $e$ is dependent on them, it
cannot be in $B$. Hence $e$ is a loop in $M_\F$, so the weight
class of $\F$ is not valid.

Conversely, assume every $F_i$ in $\F$ \emph{is} a flat of $M$.
Consider any $e \in [n]$, and find the value of $i$ such that $e
\in F_i-F_{i-1}$. After $r(F_{i-1})$ steps of the greedy
algorithm, we produce a basis of $F_{i-1}$. Since $F_{i-1}$ is a
flat, $e$ is not dependent on it, and in the next step of the
algorithm we can choose $e$. After $r$ steps, we will have a
$\w$-minimum basis of $M$ containing $e$. Therefore the weight
class of $\F$ is valid.
\end{proof}

The order complex $\Delta(L_M-\{\,\hat{0}\,, \hat{1}\,\}\,)$ is a
well understood object \cite{Bjorner}. As an immediate consequence
of Theorem \ref{theorem:main}, we get the following result.

\begin{corollary}
The Bergman complex $\B(M)$ is homotopy equivalent to a wedge of
~$\widehat{\mu}(L_M)$ $(r-2)$-dimensional spheres. Its subdivision
into weight classes is a pure, shellable simplicial complex.
\end{corollary}

Here $\widehat{\mu}(L_M) = (-1)^{r(M)} \mu_{L_M}(\hat 0, \hat 1)$
is an evaluation of the \emph{M\"{o}bius function} $\mu_{L_M}$ of
the lattice $L_M$. The M\"{o}bius function is an extremely useful
combinatorial invariant of a poset; for more information, see
\cite[Chapter 3]{Stanley}.

\medskip

\begin{figure}[h]
\centering
\includegraphics[height=2in, width=4in]{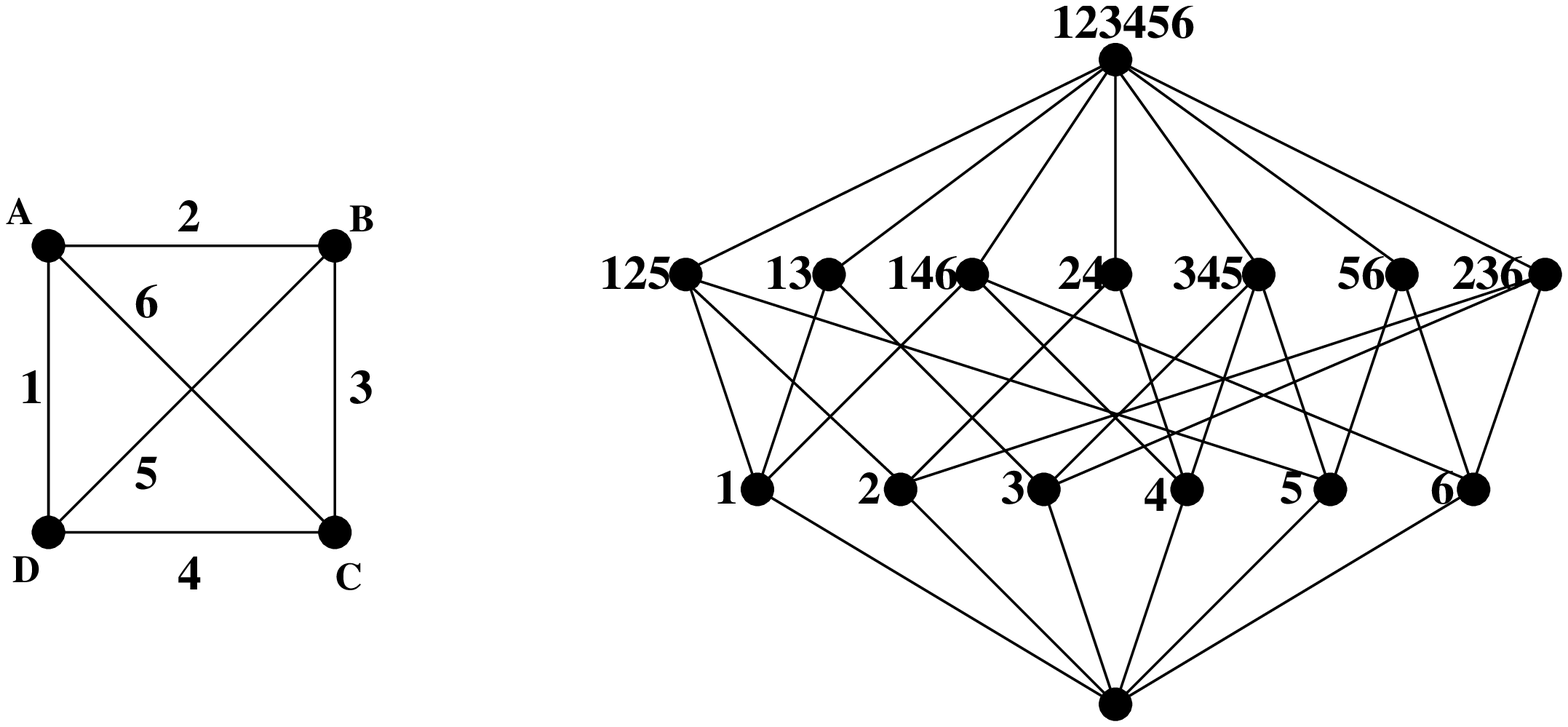}
\caption{The graph $K_4$ and its lattice of flats $\Pi_4$.}
\label{fig:K4}
\end{figure}

{\bf Example:} Let $M(K_4)$ be the graphical matroid of the
complete graph on four nodes. The bases of this matroid are given
by spanning trees. The flats are complete subgraphs and vertex
disjoint unions of complete subgraphs (see Figure~\ref{fig:K4}).
Note that in this case, the flats are in correspondence with the
partitions of the set $\{A,B,C,D\}$. In general, the flats of the
graphical matroid of $K_n$ are in bijection with partitions of the
set $[n]$.  Furthermore, the lattice of flats is the partition
lattice $\Pi_n$, which orders partitions by refinement.

\begin{figure}[h]
\centering
\includegraphics[height=2.3in ]{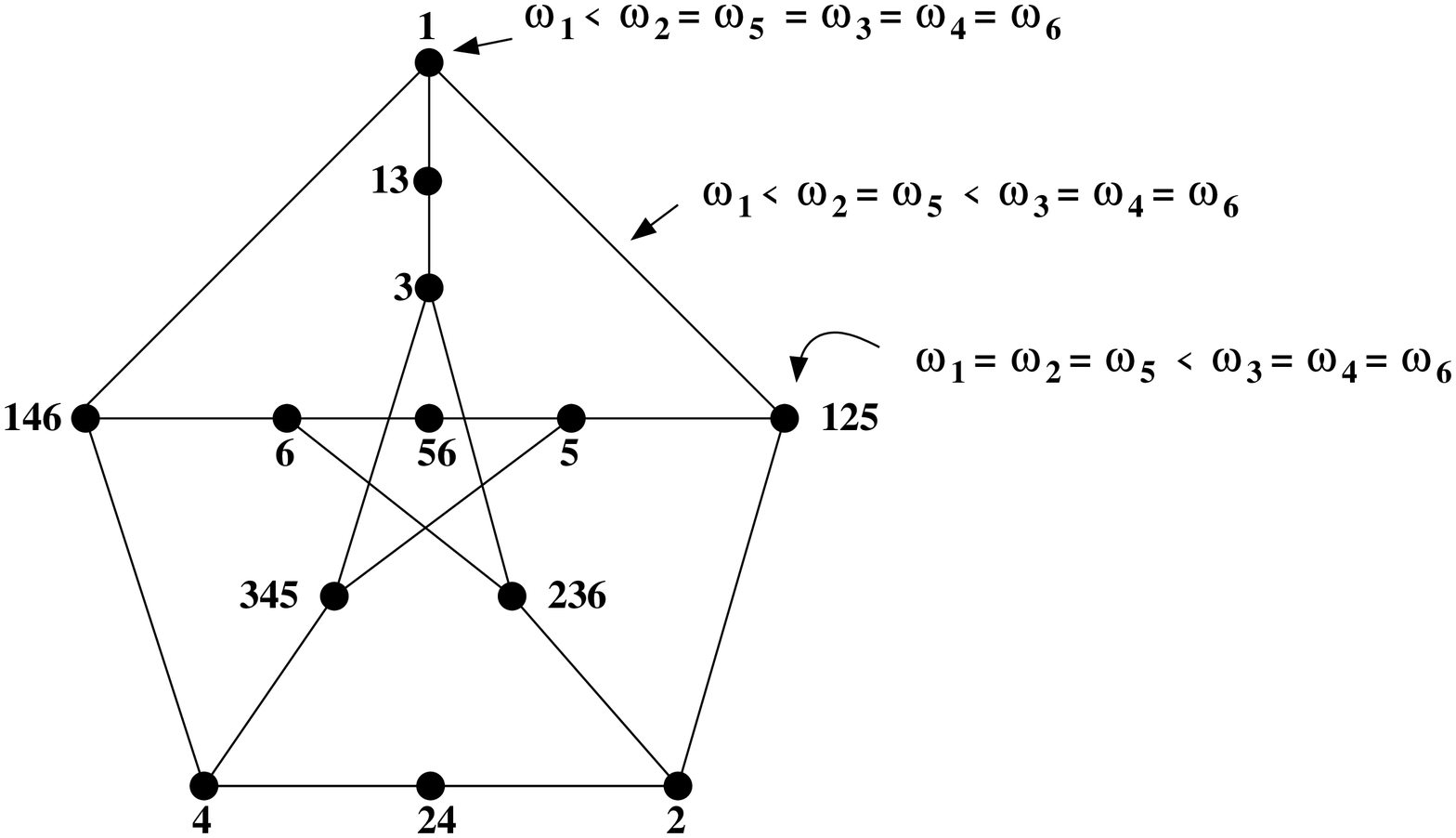}
\caption{The fine subdivision of $\B(K_4)$.} \label{fig:BK4}
\end{figure}

The fine subdivision of the Bergman complex $\B(K_4)$ is shown in
Figure~\ref{fig:BK4}. It is a wedge of six $1$-spheres. More
generally, $\B(K_n)$ is a wedge of $\widehat{\mu}(\Pi_n)=(n-1)!$
spheres of dimension $n-3$. The vertices of $\B(K_4)$ are labeled
with the corresponding flats, and a few of the corresponding
weight classes are shown. Notice that the ground set of a matroid
is always a flat, which corresponds to the weight class in which
all weights are equal. We removed this weight class when
normalizing the Bergman complex to the sphere.

The fine subdivision of the Bergman complex is almost the Petersen
graph. The only difference is the presence of the three extra
vertices, $13, 24$ and $56$. In the coarse subdivision into
$M_{\w}$-equivalence classes, these three vertices do not appear.
For example, the weight class $\omega_1 < \omega_3 < \omega_2 =
\omega_4 = \omega_5 = \omega_6$ induces the same matroid $M_{\w}$
as $\omega_1 = \omega_3 < \omega_2 = \omega_4 = \omega_5 =
\omega_6$ and $\omega_3 < \omega_1 < \omega_2 = \omega_4 =
\omega_5 = \omega_6$. Next we describe the relationship between
these two subdivisions in general.

\medskip

The coarse decomposition of $\B(M)$ into cells which induce the
same $M_{\w}$ is also a pure, polyhedral complex: it is a
subcomplex of the spherical polar to the matroid polytope of $M$.
To describe this decomposition, it is enough to describe its
full-dimensional cells.

Therefore, we only need to determine when two full-dimensional
weight classes give the same matroid $M_{\w}$. It is clearly
enough to answer this question when the two weight classes are
\emph{adjacent}; \emph{i.e.}, the intersection of their closures
is a facet of both. This happens when the two corresponding flags,
which have one flat in each rank, are equal in all but one rank.

Let $\diamondsuit$ be the \emph{diamond poset}; \emph{i.e.}, the
rank $2$ poset consisting of a minimum element, a maximum element,
and two rank $1$ elements.

\begin{ntheorem}\label{theorem:subdiv}
Suppose that the weight classes of two maximal flags $\F$ are $\F'$
are adjacent. Say $\F$ and $\F'$ only differ in rank $i$; that is,
$\F-F_i = \F' - F_i'$. Then the following conditions are equivalent:

\begin{enumerate}
\item[(i)] $M_{\F} = M_{\F'}$. \item[(ii)] $M_{\F} = M_{\F-F_i}$.
\item[(iii)] $F_i \cup F_i' = F_{i+1}$. \item[(iv)] The interval
$[F_{i-1},F_{i+1}]$ of $L_M$ is a diamond poset.
\end{enumerate}
\end{ntheorem}

\begin{proof}
Let $M_j = (M|\,F_j) / F_{j-1}$, $M_j' = (M|\,F_j') / F_{j-1}'$,
$N_i=(M|F_{i+1})/F_{i-1}$, and $N = M_1 \oplus \cdots \oplus
M_{i-1} \oplus M_{i+2} \oplus \cdots \oplus M_{k+1}$. By
Proposition \ref{prop:Mw2},
\[
M_{\F} = N  \oplus M_i \oplus M_{i+1}, \quad M_{\F'} = N \oplus
M_i' \oplus M_{i+1}',  \quad M_{\F-F_i} = N \oplus N_i.
\]
Since $M_i, M_{i+1}, M_i'$ and $M_{i+1}$ have rank $1$ and $N_i$
has rank $2$,
\begin{eqnarray*}
L_{M_i \oplus M_{i+1}} &=& \{\emptyset, F_i-F_{i-1}, F_{i+1}-F_i,
F_{i+1}-F_{i-1}\} \cong \diamondsuit, \\
 L_{M_i' \oplus M_{i+1}'} &=& \{\emptyset,
F_i'-F_{i-1}, F_{i+1}-F_i', F_{i+1}-F_{i-1}\} \cong \diamondsuit, \\
L_{N_i} &=& \{F-F_{i-1} : F \in [F_{i-1},F_{i+1}]\} \cong
[F_{i-1},F_{i+1}].
\end{eqnarray*}


If $(iv)$ does not hold, then we know immediately that $L_{N_i}
\neq L_{M_i \oplus M_{i+1}}$. Also $F_i \cup F_i' \neq F_{i+1}$,
and therefore $L_{M_i \oplus M_{i+1}} \neq L_{M_i' \oplus
M_{i+1}'}$.

If $(iv)$ holds, then $F_i$ and $F_i'$ are the only rank $i$ flats
of $M$ in $[F_{i-1},F_{i+1}]$. Since $N_i$ has no loops, $(iii)$ holds;
and therefore $L_{M_i \oplus M_{i+1}} = L_{M_i' \oplus
M_{i+1}'}=L_{N_i}$.
\end{proof}

\section{The space of phylogenetic trees} \label{section:trees}

In this section, we show that the Bergman fan $\widetilde{\B}(K_n)$ of the
matroid of the
complete graph $K_n$ is homeomorphic to the space of phylogenetic
trees $\T_n$, as defined in \cite{Billera}. To do so, we start by
reviewing the connection between phylogenetic trees and ultrametrics.

\begin{definition}
A \emph{dissimilarity map} on $[n]$ is a map $\delta: [n] \times
[n] \rightarrow \R$ such that $\delta(i,i) = 0$ for all $i \in
[n]$, and $\delta(i,j) = \delta(j,i)$ for all $i,j \in [n]$.
\end{definition}

\begin{definition}
A dissimilarity map is an \emph{ultrametric} if, for all $i,j,k
\in [n]$, two of the values $\delta(i,j), \delta(j,k)$ and
$\delta(i,k)$ are equal and not less than the third.
\end{definition}

\begin{figure}[h]
\centering
\includegraphics[height=1.3in ]{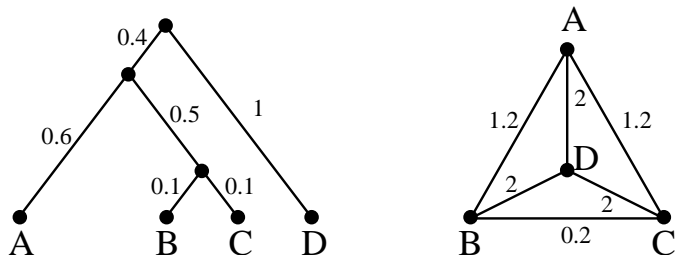}
\caption{An equidistant tree and its corresponding distance
function.} \label{figure:equitree}
\end{figure}

An \emph{equidistant $n$-tree} $T$ is a rooted tree with $n$
leaves labelled $1, \ldots, n$, and lengths assigned to each edge
in such a way that the path from the root to any leaf has the same
length. The internal edges are forced to have positive lengths,
while the edges incident to a leaf are allowed to have negative
lengths. Figure \ref{figure:equitree} shows an example of an
equidistant $4$-tree.

To each equidistant $n$-tree $T$ we assign a distance function
$d_T :[n] \times [n] \rightarrow \R$: the distance $d_T(i,j)$ is
equal to the length of the path joining leaves $i$ and $j$ in $T$.
Such a distance function can be regarded as a weight vector on the
edges of $K_4$; Figure \ref{figure:equitree} also shows the
distance function of the tree shown.


We can think of equidistant trees as a model for the evolutionary
relationships between a certain set of species. The various
species, represented by the leaves, descend from a single root.
The descent from the root to a leaf tells us the history of how a
particular species branched off from the others until the present
day. For more information on the applications of this and other
similar models, see for example \cite{Billera} and \cite{Semple}.

The connection between equidistant trees and ultrametrics is given
by the following theorem.

\begin{theorem} \cite[Theorem 7.2.5]{Semple}
A map $\delta: [n] \times [n] \rightarrow \R$ is an ultrametric if
and only if it is the distance function of an equidistant
$n$-tree.
\end{theorem}


We can now explain the relationship between the Bergman fan
$\widetilde{\B}(K_n)$ and phylogenetic trees.

\begin{ntheorem}\label{theorem:trees}
A dissimilarity map $\delta: [n] \times [n] \rightarrow \R$ is an
ultrametric if and only if the corresponding weight vector on the
edges of $K_n$ is in the Bergman fan $\widetilde{\B}(K_n)$.
\end{ntheorem}

\begin{proof}
We claim that the following three statements about a weight function on the
edges of $K_n$ are equivalent.
\begin{enumerate}
\item[$(i)$] In any triangle, the largest weight is achieved (at
least) twice. \item[$(ii)$] In any cycle, the largest weight is
achieved (at least) twice. \item[$(iii)$] Every edge is in a
spanning tree of minimum weight.
\end{enumerate}

The theorem will follow from this claim, because ultrametrics are
characterized by $(i)$ and weight functions in the Bergman complex
are characterized by $(iii)$.

The implication $(ii) \Rightarrow (i)$ is trivial. Conversely,
assume that $(i)$ holds and $(ii)$ does not. Without loss of
generality, assume that the cycle $v_1 v_2 \ldots v_k$ has
$v_1v_2$ as its unique edge of largest weight. The largest weight
in triangle $v_1v_2v_3$ must be achieved at $\w(v_1v_2) =
\w(v_1v_3)$. The largest weight in triangle $v_1v_3v_4$ must then
be achieved at $\w(v_1v_3) = \w(v_1v_4)$. Continuing in this way
we get that $\w(v_1v_2) = \w(v_1v_3) = \cdots = \w(v_1v_k)$, and
$(ii)$ follows.

Now we prove $(ii) \Rightarrow (iii)$. Consider an arbitrary edge
$f$. Let $T$ be a spanning tree of minimum weight. If $f \in T$ we
are done; otherwise, $T \cup f$ has a unique cycle. There is at
least one edge $e$ in this cycle with $\w(e) \geq \w(f)$.
Therefore, the weight of the spanning tree $T \backslash e \cup f$
is not larger than the weight of $T$. This is then a spanning tree
of minimum weight containing $f$.

Finally, assume that $(iii)$ holds and $(i)$ does not. Assume that
the triangle with edges $e,f,g$ has $\w(e) > \w(f), \w(g)$, and
consider a spanning tree $T$ of minimum weight which contains edge
$e$. If $f$ is in $T$, then $g$ cannot be in $T$, and replacing
$e$ with $g$ will give a spanning tree of smaller weight.  If
neither $f$ nor $g$ is in $T$, we can still replace $e$ with one
of them to obtain a spanning tree of smaller weight. If we could
not, that would imply that both $f$ and $g$ form a cycle when
added to $T \backslash e$. Call these cycles $C_f$ and $C_g$. But
then $(C_f \backslash f) \cup (C_g \backslash g) \cup e$ would
contain a cycle in $T$, a contradiction.
\end{proof}

The previous two theorems give us a one-to-one correspondence
between the vectors in the Bergman fan $\widetilde{\B}(K_n)$ and the
equidistant $n$-trees: $\widetilde{\B}(K_n)$ parameterizes
equidistant $n$-trees by the distances between their leaves. This
leads us to consider the space of trees
$\T_n$ of~\cite{Billera}. This space parameterizes equidistant
$n$-trees in a different way: it keeps track of their
combinatorial type, and the lengths of their internal edges. We
recall the construction of the space $\T_n$.  Each maximal cell
corresponds to a combinatorial type of rooted binary tree on $n$
labeled leaves; \emph{i.e.}, a rooted tree where each internal
vertex has two descendants.  Such trees have $n-2$ internal edges,
and are parameterized by vectors in $\R_{>0}^{\,\,n-2}$ recording
these edge lengths. Moving to a lower dimensional face of a maximal
cell corresponds to setting some of these edge lengths to $0$,
which gives non-binary degenerate cases of the original tree.
Maximal cells are glued along these lower-dimensional cells when
two trees specialize to the same degenerate tree.

\begin{figure}[h]
\centering
\includegraphics[height=.9in]{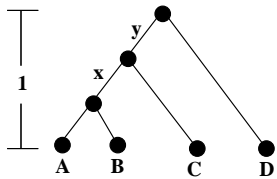}
\caption{Combinatorial type of tree with 4 leaves.}
\label{fig:combo}
\end{figure}

Given a fixed combinatorial type of tree and the vector of
internal edge lengths, we can recover the pairwise distances of
leaves as linear functions on the internal edge lengths.  For
example, consider the tree type of Figure~\ref{fig:combo}.  We
obtain
$(\delta(A,B),\delta(A,C),\delta(A,D),\delta(B,C),\delta(B,D),\delta(C,D))
\in \B(K_4)$ from $(x,y)$ by the map $f$ : $(x,y)$ $\mapsto$
$(2(1-x-y),2(1-y),2,2(1-y),2,2)$.  The converse is also true;
given the pairwise distances of leaves we can recover the internal
edge lengths via linear relations on these
distances~\cite{Semple}.

In general, doing this for each type of tree, we get a map $f:\T_n
 \rightarrow \widetilde{\B}(K_n)$.  It follows from the previous
 two theorems that $f$ is a one-to-one correspondence between $\T_n$ and
 $\widetilde{\B}(K_n)$. We will now see that, in fact, $\T_n$ and
 $\widetilde{\B}(K_n)$ have the same combinatorial structure.

\begin{proposition}
The map $f:\T_n \rightarrow \widetilde{\B}(K_n)$ is a piecewise
linear homeomorphism. It identifies the decomposition of the space
of trees $\T_n$ into combinatorial tree types with the coarse
decomposition of the Bergman fan $\widetilde{\B}(K_n)$.
\end{proposition}

\begin{proof}
Restricting to a maximal cell of $\T_n$, corresponding to a fixed
tree type, $f$ is a linear map from the lengths of internal edges
(in the space of trees) to the pairwise distances of the leaves
(in the space of ultrametrics). Also, it is clear that when two
maximal cells of $\T_n$ intersect, the linear restrictions of $f$
to these two cells agree on their intersection. The first claim
follows.

Suppose we are given a combinatorial type of equidistant $n$-tree.
From the branching order of each triple of leaves (\emph{i.e.},
which, if any, of the three branched off first), we can recover
which edges of each triangle of $K_n$ are maximum in the
corresponding weight vector. In turn, this allows us to recover
which edges of any cycle are maximum: one can check that an edge
is maximum in a cycle $C$ if and only if it is maximum in each
triangle that it forms with a vertex of $C$.  Knowing the maximum
edges of each cycle of the graph, we can determine $M_{\omega}$
using the following version of the greedy algorithm. Start with
the complete graph $K_n$ and break its cycles successively: at
each step pick an existing cycle, and remove one of its maximum
edges. The trees which can result by applying this procedure are
precisely the $\w$-minimum spanning trees~\cite{Kozen}. Therefore
$f$ maps a fixed tree type class of $\T_n$ to a fixed
$M_{\w}$-equivalence class; \emph{i.e.}, a fixed cell in the coarse
subdivision of $\widetilde{\B}(K_n)$.

Conversely, suppose we are given $M_{\omega}$ (which has no loops)
and we want to determine the combinatorial tree type of
$f^{-1}(\w)$. Consider the edges $\{e,f,g\}$ of any
triangle in $K_n$; we can find out whether $e$ is maximum in this
triangle as follows. Take a minimum spanning tree $T$ containing
$e$. Either $T \backslash e \cup f$ or $T\backslash e \cup g$ is a
spanning tree; assume it is the first. If $T \backslash e \cup f$
is a minimum spanning tree, then $\omega(e) = \omega(f)$, and $e$
is maximum in the triangle. Otherwise $\omega(e) < \omega(f)$
and $e$ is not maximum in the triangle. Determining this
information for each triangle tells us, for each triple of leaves,
which one (if any) branched off first in the corresponding tree.
It is easy to reconstruct the combinatorial type of the tree from
this data, in the same way that one recovers an equidistant tree
from its corresponding ultrametric \cite[Theorem 7.2.5]{Semple}.
\end{proof}

The link of the origin in the coarse subdivision of $\T_n$, which
we call $T_n$, is a simplicial complex which has appeared in many
different contexts. It was first considered by Boardman
\cite{Boardman}, and also studied by Readdy \cite{Readdy},
Robinson and Whitehouse \cite{Robinson}, Sundaram \cite{Sundaram},
Trappmann and Ziegler \cite{Trappmann}, Vogtmann \cite{Vogtmann},
and Wachs \cite{Wachs}, among others. By Theorem
\ref{theorem:main}, the link of the origin in the fine subdivision
of $\widetilde{\B}(K_n)$ is the order complex of the proper part
of the partition lattice $\Pi_n$. We conclude the following
result.

\begin{corollary}
The order complex of the proper part of the partition lattice
$\Pi_n$ is a subdivision of the complex $T_n$.
\end{corollary}

This provides a new explanation of the known result ~\cite{Readdy,
Robinson,Sundaram, Trappmann, Vogtmann,Wachs} that these two
simplicial complexes are homotopy equivalent; namely, they have
the homotopy type of a wedge of $(n-1)!$ $(n-3)$-dimensional
spheres.


\begin{figure}[H]
\centering \includegraphics[height=2.8in]{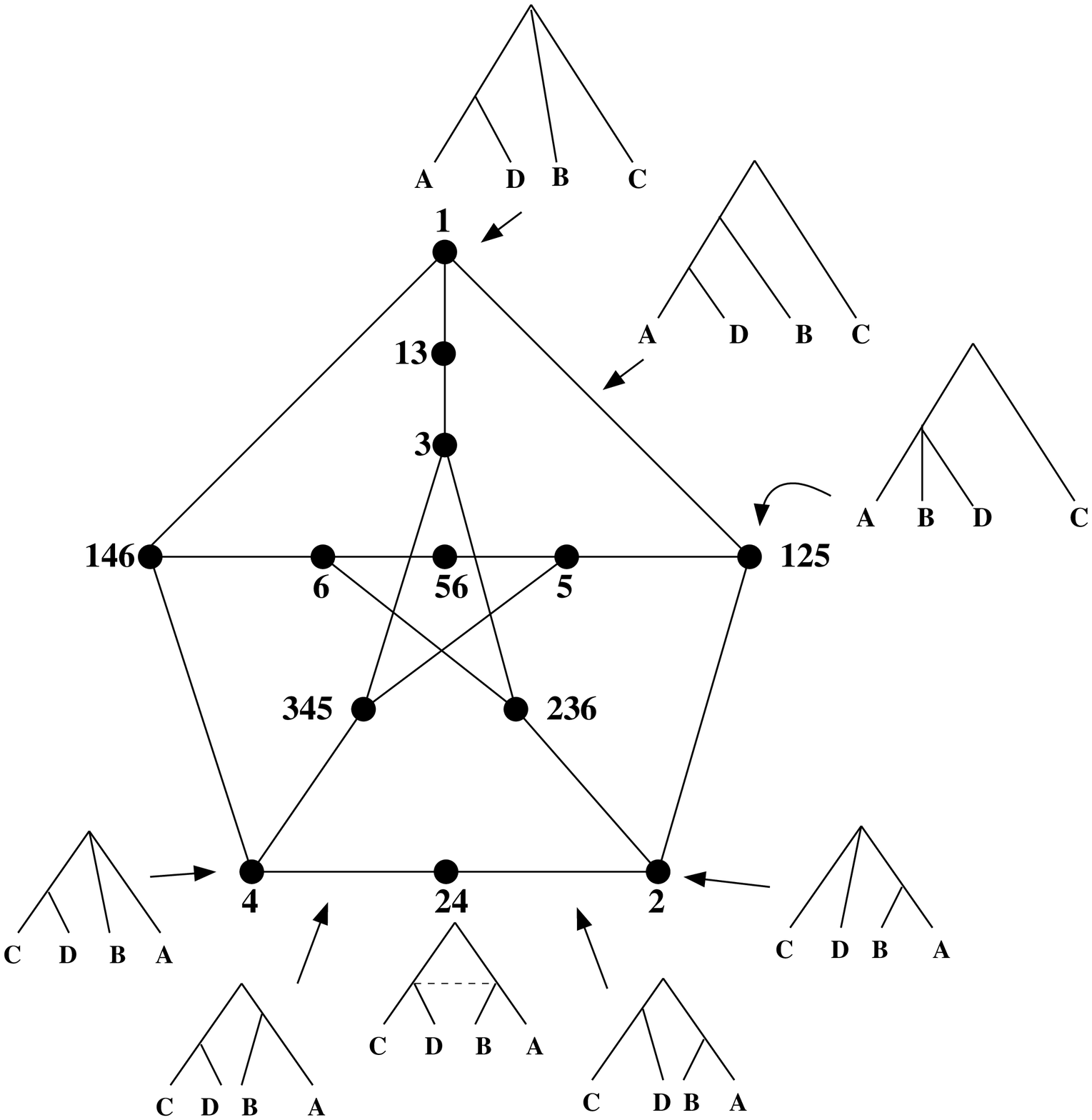} \caption{The fine
subdivision of $\B(K_4)$ revisited.} \label{fig:TK4}
\end{figure}

\medskip

Let us now revisit the example of the last section. In
Figure~\ref{fig:TK4} we show the Bergman complex $\B(K_4)$, with
some of the corresponding trees. We now know that this is a
subdivision of $T_4$, the link of the origin in the space of
phylogenetic trees with $4$ leaves, which is the Petersen graph.
The three extra vertices in the fine subdivision are $13$, $24$
and $56$. The tree corresponding to vertex $24$ of the fine
subdivision has the property that the vertex joining the leaves
$C$ and $D$ is at the same height as the vertex joining the leaves
$A$ and $B$. This information is not captured by the combinatorial
type of the tree; \emph{i.e.}, by the coarse subdivision.


\begin{figure}[H]
\centering
\includegraphics[height=3in, width=4in]{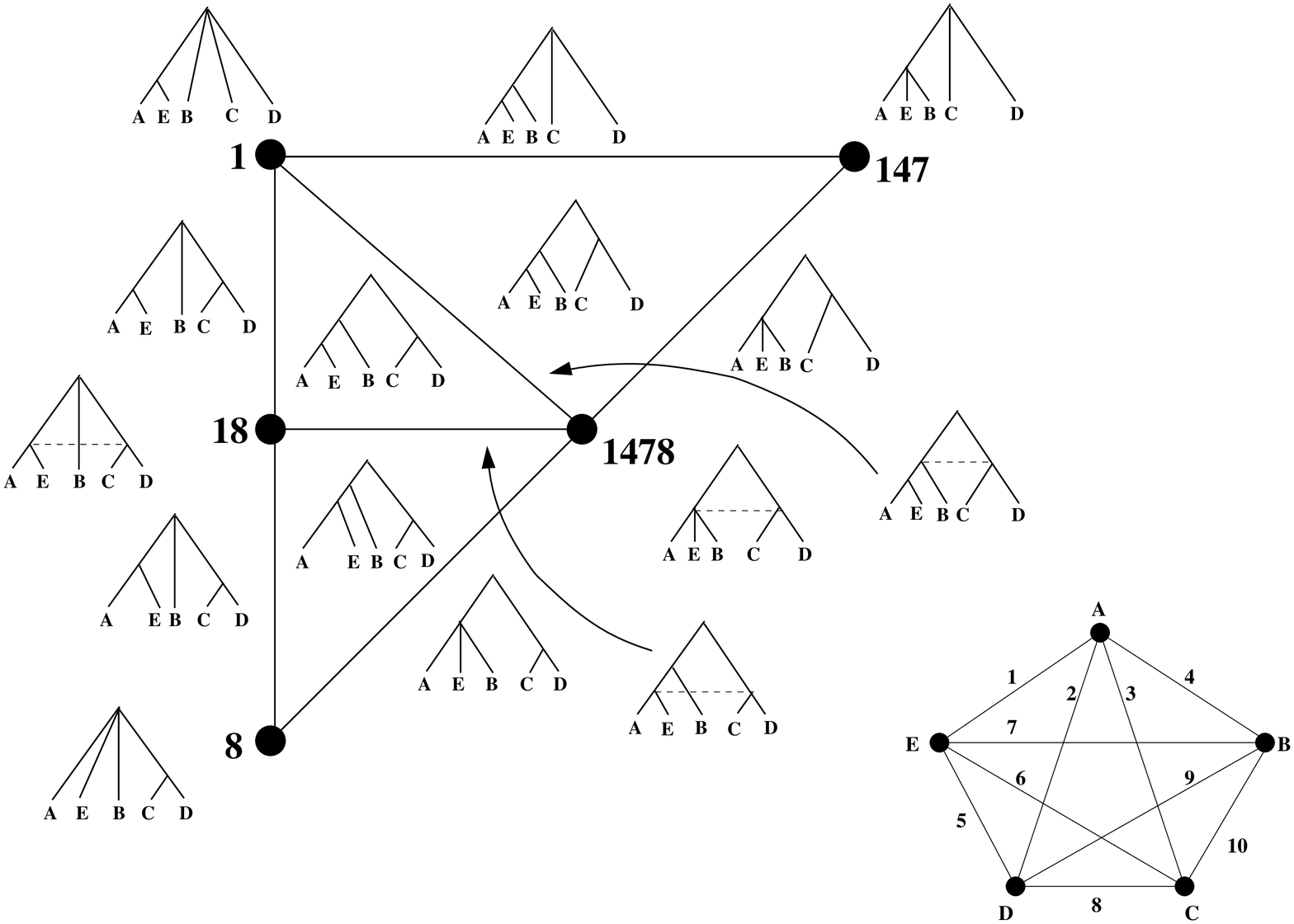}
\caption{A piece of the fine subdivision of $\B(K_5)$.}
\label{fig:NBK5}
\end{figure}

In Figure~\ref{fig:NBK5}, we show a representative piece of the
fine subdivision of the space of trees with $5$ leaves, with $K_5$
labeled as shown.

\section{Acknowledgments}

We would like to thank Bernd Sturmfels for suggesting to us the
project of studying the Bergman complex of a matroid, and for many
helpful conversations on this topic. We would also like to thank
Lou Billera, G\"{u}nter Ziegler, and an anonymous referee, for
reading preliminary versions of this manuscript and helping us
improve the presentation.

\medskip

\textsc{Mathematical Sciences Research Institute, Berkeley, CA,
USA.}


\small

\begin{thebibliography}{99}

\bibitem{Bergman}
G. Bergman.  \emph{The logarithmic limit-set of an algebraic
variety.}  Trans. Amer. Math. Soc. {\bf 157} (1971) 459-469.

\bibitem{Bieri}
R. Bieri and J. Groves. \emph{ The geometry of the set of
characters induced by valuations.} J. Reine Angew. Math. {\bf 347}
(1994) 168-195.

\bibitem{Billera}
L. Billera, S. Holmes, and K. Vogtmann.  \emph{Geometry of the
space of phylogenetic trees.}  Adv. in Appl. Math. {\bf 27} (2001)
733-767.

\bibitem{Bjorner}
A. Bj\"{o}rner. \emph{The homology and shellability of matroids
and geometric lattices.} Matroid Applications.  Encyclopedia of
Mathematics and Its Applications, vol. 40. Cambridge University
Press, Cambridge, 1992.

\bibitem{Boardman}
J. M. Boardman. \emph{Homotopy structures and the language of
trees.} Proc. Sympos. Pure Math. {\bf 22} (1971), 37-58.

\bibitem{Borovik}
A. Borovik, I. Gelfand, and N. White. \emph{Coxeter matroids.}
Birkhauser, Boston, 2003.

\bibitem{Einsiedler}
M. Einsiedler, M. Kapranov, and D. Lind. \emph{Non-archimedean
amoebas and tropical varieties,} preprint,
\textsf{arXiv:math.AG/0408311}, 2004.

\bibitem{Kozen}
D. Kozen. \emph{The design and analysis of algorithms.} Lecture
Notes in Computer Science, Springer-Verlag, 1991.

\bibitem{Oxley}
J. G. Oxley. \emph{Matroid theory.} Oxford University Press, New
York, 1992.

\bibitem{Readdy}
M. A. Readdy. \emph{The pre-WDVV ring of physics and its
topology,} to appear in Ramanujan J., 2005.

\bibitem{Robinson}
A. Robinson and S. Whitehouse. \emph{The tree representation of
$\sigma_{n+1}$.}  J. Pure Appl. Algebra {\bf 111} (1996) 245-253.

\bibitem{Semple}
C. Semple and M. Steel. \emph{Phylogenetics.} Oxford Lecture
Series in Mathematics and Its Applications, Oxford University
Press, 2003.

\bibitem{Speyer}
D. Speyer and B. Sturmfels. \emph{The tropical Grassmannian.} Adv.
in Geom. {\bf 4} (2004), 389-411.

\bibitem{Stanley}
R. P. Stanley. \emph{Enumerative combinatorics, vol.1.} Cambridge
University Press, New York, 1986.

\bibitem{Sturmfels}
B. Sturmfels. \emph{Solving systems of polynomial equations.} CBMS
Regional Conference Series in Mathematics, vol. 97.  American
Mathematical Society, Providence, RI, 2002.

\bibitem{Sundaram}
S. Sundaram. \emph{Homotopy of the non-modular partitions and the
Whitehouse module}.  J. Algebraic Combin. {\bf 9} (1999) 251--269.

\bibitem{Trappmann}
H. Trappmann and G. M. Ziegler.  \emph{Shellability of complexes
of trees.} J. Combin. Theory Ser. A {\bf 82} (1998) 168-178.

\bibitem{Vogtmann}
K. Vogtmann.  \emph{Local structure of some $OUT(F_n)$-complexes.}
Proc. Edinb. Math. Soc. (2) {\bf 33} (1990) 367-379.

\bibitem{Wachs}
M. Wachs, personal communication, 2003.


\end{thebibliography}
\end{document}